# Study about a Differential Equation in an Infinite Servers Queue System with Poisson Arrivals Busy Cycle Distribution Study

Manuel Alberto M. Ferreira

**Abstract**

In $M|G|\infty$ queue real life practical applications, the busy period and the busy cycle probabilistic study is of main importance. But it is a very difficult task. In this chapter, we show that by solving a Riccati equation induced by this queue transient probabilities monotony study as time functions, we obtain a collection of service length distribution functions, for which both the busy period and the busy cycle have lengths with quite simple distributions, generally given in terms of exponential distributions and the degenerate at the origin distribution.

**Introduction**

In an $M|G|\infty$ queue system, customers arrive according to a Poisson process at rate $\lambda$, receive a service which length is a positive random variable with distribution function $G(\cdot)$ and mean $\alpha = \int_0^\infty [1 - G(v)]dv$. There are always available servers. The traffic intensity is $\rho = \lambda\alpha$.

In the $M|G|\infty$ queue system, such as in any other queue system, there is a sequence of idle periods and busy periods. An idle period followed by a busy period is called busy cycle.

Be, for the $M|G|\infty$ queue system, $N(t)$ the number of occupied servers or, equivalently, the number of customers being served at time t. According to [1], $p_{00}(t) = P(N(t) = 0|N(0) = 0)$ is given by

$$p_{00}(t) = e^{-\lambda \int_0^t [1-G(v)]dv} \quad (1).$$

Be also $p_{1'0}(t) = P(N(t) = 0|N(0 - \varepsilon) = 0), \forall_{\epsilon>0}$ and $N(0) \geq 1$, meaning $N(0 - \varepsilon) = 0, \forall_{\epsilon>0}$ and $N(0) \geq 1$ that the time origin is an instant at which a busy period begins. Considering two subsystems, one corresponding to the initial customer, and the other to the other infinite servers, obviously behaving in an independent way, we conclude that

$$p_{1'0}(t) = p_{00}(t)G(t) \quad (2).$$

As $\frac{d}{dt}p_{1'0}(t) = p_{00}(t)(1 - G(t))\left(\frac{g(t)}{1-G(t)} - \lambda G(t)\right)$, where $g(t) = \frac{d}{dt}G(t)$, we conclude that the monotony of $p_{1'0}(t)$, as time function, depends on the sign of the expression

$$\frac{g(t)}{1-G(t)} - \lambda G(t) \quad (3).$$

We will see in this paper, how (3) induces a Riccati equation, which solution, under specific restrictions, leads to the determination of a collection of distributions of the service length for which either the busy period or the busy cycle, that is: the respective lengths are governed for simple probability laws.

So, we intend to solve that equation, obtaining a collection of the service length distribution functions, for which we will determine the respective $M|G|\infty$ busy period length and busy cycle distribution functions. Here plays a fundamental role the busy period length Laplace transform, see [2].

This chapter is the enlarged and updated version of the paper published in the peer-reviewed proceedings intitled, in Portuguese, *Estatística Jubilar-Actas do XII Congresso Anual da Sociedade Portuguesa de Estatística,* see [3].

**Busy Period and Busy Cycle Probabilistic Study**

Defining $\beta(t) = \frac{g(t)}{1-G(t)} - \lambda G(t)$, after some operations, it can be given the form

$$\frac{dG(t)}{dt} = -\lambda G^2(t) - (\beta(t) - \lambda)G(t) + \beta(t) \quad (4),$$

that is the one of a Riccati equation which solution is

$$G(t) = 1 - \frac{1}{\lambda}\frac{(1-e^{-\rho})e^{-\lambda t-\int_0^t \beta(u)du}}{\int_0^\infty e^{-\lambda w-\int_0^w \beta(u)du}dw - (1-e^{-\rho})\int_0^t e^{-\lambda w-\int_0^w \beta(u)du}dw}, t \geq 0, -\lambda \leq \frac{\int_0^t \beta(u)du}{t} \leq \frac{\lambda}{e^\rho - 1} \quad (5).$$

Note that $G(t) = 1, t \geq 0$ is a solution contemplated in (5) for $\frac{\int_0^t \beta(u)du}{t} = -\lambda$. After (5) and

$$\bar{B}(s) = 1 + \lambda^{-1}\left(s - \frac{1}{\int_0^\infty e^{-st-\lambda\int_0^t[1-G(v)]dv}dt}\right) \quad (6),$$

that is the $M|G|\infty$ system busy period length Laplace transform, see [2], we obtain

$$\bar{B}(s) = \frac{1 - (s+\lambda)(1-G(0))L\left(e^{-\lambda t - \int_0^t \beta(u)du}\right)}{1 - \lambda(1-G(0))L\left(e^{-\lambda t - \int_0^t \beta(u)du}\right)}, -\lambda \le \frac{\int_0^t \beta(u)du}{t}$$

$$\le \frac{\lambda}{e^\rho - 1}, \text{ where } L(\cdot) \text{ means Laplace transform and}$$

$$G(0) = \frac{\lambda \int_0^\infty e^{-\lambda w - \int_0^w \beta(u)du} dw + e^{-\rho} - 1}{\lambda \int_0^\infty e^{-\lambda w - \int_0^w \beta(u)du} dw} \quad (7).$$

After (7) we can compute $\frac{1}{s}\bar{B}(s)$ which inversion leads to

$$B(t) = (1 - (1-G(0))(e^{-\lambda t - \int_0^t \beta(u)du} + \lambda \int_0^t e^{-\lambda w - \int_0^w \beta(u)du} dw)) * \sum_{n=0}^\infty \lambda^n (1-G(0))^n \left(e^{-\lambda t - \int_0^t \beta(u)du}\right)^{*n}, t \ge 0, -\lambda \le \frac{\int_0^t \beta(u)du}{t} \le \frac{\lambda}{e^\rho - 1} \quad (8),$$

the $M|G|\infty$ queue system busy period length distribution function, being $*$ the convolution operator.

If $\beta(t) = \beta$ (constant), (5) becomes

$$G(t) = 1 - \frac{(1-e^{-\rho})(\lambda+\beta)}{\lambda e^{-\rho}(e^{(\lambda+\beta)t}-1)+\lambda}, t \ge 0, -\lambda \le \beta \le \frac{\lambda}{e^\rho - 1} \quad (9)$$

and (8),

$$B^{(\beta)}(t) = 1 - \frac{\lambda+\beta}{\lambda}(1-e^{-\rho})e^{-e^{-\rho}(\lambda+\beta)t}, t \ge 0, -\lambda \le \beta \le \frac{\lambda}{e^\rho - 1} \quad (10).$$

So, if the service length distribution is given by (9), the busy period length distribution is the mixture of a degenerate at the origin distribution with an exponential distribution.

For $\beta = -\lambda$,

$$G(t) = 1, t \geq 0 \text{ and } B^{(-\lambda)}(t) = 1, t \geq 0,$$

that is: if the service length distribution is degenerate at the origin with probability 1, the same happens with the busy period length distribution.

Note also that if $\beta = \frac{\lambda}{e^\rho - 1}$,

$$B^{(\frac{\lambda}{e^\rho-1})} = 1 - e^{-\frac{\lambda}{e^\rho-1}t}, t \geq 0$$

(purely exponential). And $B(t)$, given by (8) satisfies

$$B(t) \geq 1 - e^{-\frac{\lambda}{e^\rho-1}t}, t \geq 0, -\lambda \leq \frac{\int_0^t \beta(u)du}{t} \leq \frac{\lambda}{e^\rho-1}.$$

Calling the busy cycle length Laplace transform $\bar{Z}(s)$, knowing that for the $M|G|\infty$ queue the idle period length is exponentially distributed with parameter $\lambda$, as it happens with any queue system with Poisson process arrivals at rate $\lambda$, and that in the $M|G|\infty$ queue the idle period and the busy period are independent, see [4],

$$\bar{Z}(s) = \frac{\lambda}{\lambda+s}\bar{B}(s) \quad (11).$$

After (7) and (11), computing $\frac{1}{s}\bar{Z}(s)$, and then inverting, we get:

$$Z(t) = (\lambda e^{-\lambda t}) * (1 - (1 - G(0))(e^{-\lambda t - \int_0^t \beta(u)du} + \lambda \int_0^t e^{-\lambda w - \int_0^w \beta(u)du} dw)) *$$
$$\sum_{n=0}^{\infty} \lambda^n (1 - G(0))^n \left(e^{-\lambda t - \int_0^t \beta(u)du}\right)^{*n}, t \geq 0, -\lambda \leq \frac{\int_0^t \beta(u)du}{t} \leq \frac{\lambda}{e^\rho-1} \quad (12),$$

the $M|G|\infty$ system busy cycle length distribution function.

If $\beta(t) = \beta$ (constant), (12) becomes

$$Z^{(\beta)}(t) = 1 - \frac{(1-e^{-\rho})(\lambda+\beta)}{\lambda - e^{-\rho}(\lambda+\beta)} e^{-e^{-\rho}(\lambda+\beta)t} + \frac{\beta}{\lambda - e^{-\rho}(\lambda+\beta)} e^{-\lambda t}, t \geq 0, -\lambda \leq \beta \leq \frac{\lambda}{e^\rho-1} \quad (13).$$

So, if the service length distribution is given by (9), the busy cycle length distribution is the mixture of two exponential distributions.

For $\beta = -\lambda$,

$$G(t) = 1, t \geq 0 \text{ and } Z^{(-\lambda)}(t) = 1 - e^{-\lambda t}, t \geq 0.$$

that is: if the service length distribution is degenerate at the origin with probability 1, the busy cycle length distribution is exponential with parameter $\lambda$: the same as the idle period, as expected.

And $Z(t)$, defined by (12), satisfies

$$Z(t) \leq 1 - e^{-\lambda t}, t \geq 0, -\lambda \leq \frac{\int_0^t \beta(u)du}{t} \leq \frac{\lambda}{e^\rho - 1}.$$

Note also that if $\beta = \frac{\lambda}{e^\rho - 1}$,

$$Z^{(\frac{\lambda}{e^\rho - 1})}(t) = 1 - \frac{(e^\rho - 1)e^{-\frac{\lambda}{e^\rho - 1}t} - e^{-\lambda t}}{e^\rho - 2}, t \geq 0.$$

This formula is coherent even for $\rho = \ln 2$, since $\lim_{\rho \to \ln 2}\left(1 - \frac{(e^\rho - 1)e^{-\frac{\lambda}{e^\rho - 1}t} - e^{-\lambda t}}{e^\rho - 2}\right) =$

$\lim_{\rho \to \ln 2} \frac{e^\rho - 2 - (e^\rho - 1)e^{-\frac{\lambda}{e^\rho - 1}t} + e^{-\lambda t}}{e^\rho - 2} = 1 - (1 + \lambda t)e^{-\lambda t}.$

And $Z(t)$ given by (12) satisfies

$$Z(t) \geq 1 - \frac{(e^\rho - 1)e^{-\frac{\lambda}{e^\rho - 1}t} - e^{-\lambda t}}{e^\rho - 2}, t \geq 0, -\lambda \leq \frac{\int_0^t \beta(u)du}{t} \leq \frac{\lambda}{e^\rho - 1}.$$

For more information about the $M|G|\infty$ queue system busy period length and busy cycle length distribution functions circumstances of exponential behavior see, for instance, [5].

## Conclusions

The $p_{1'0}(t)$ monotony study, as time function, for the $M|G|\infty$ queue system, leads to the consideration of a Riccati equation. Its solution, under specific restrictions, is a collection of service length distributions for which either the busy period or the busy cycle length distributions have quite simple mathematical expressions, generally given in terms of exponential distributions and the degenerate at the origin distribution, which is extremely rare in queues occupation study.

## Acknowledgements

This work is financed by national funds through FCT - Fundação para a Ciência e Tecnologia, I.P., under the project FCT UIDB/04466/2022. Furthermore, the author thanks the ISCTE-IUL and ISTAR-IUL, for their support.